\documentclass[12pt,twoside]{article}
\usepackage{amsmath, amssymb, amsfonts, amstext}

\date{}

\begin{document}

\centerline{\Large{\bf Extremal decomposition problems in the Euclidean space}}

\centerline{}

\centerline{\bf {K.A. Gulyaeva}\footnote{Far Eastern Federal University, Vladivostok, Russia, kgulyayeva@gmail.com}, S.I. Kalmykov\footnote{Far Eastern Federal University, Vladivostok, Russia, Institute of
Applied Mathematics, Far Eastern Branch of the Russian Academy of Sciences, Vladivostok, Russia, sergeykalmykov@inbox.ru}, E.G. Prilepkina\footnote{Far Eastern Federal University, Vladivostok, Russia, Institute of Applied Mathematics, Far Eastern Branch of the Russian Academy of
Sciences, Vladivostok, Russia, pril-elena@yandex.ru}}

\newtheorem{Theorem}{\quad Theorem}[section]

\newtheorem{Definition}[Theorem]{\quad Definition}

\newtheorem{Corollary}[Theorem]{\quad Corollary}

\newtheorem{Lemma}[Theorem]{\quad Lemma}

\newtheorem{Example}[Theorem]{\quad Example}

\centerline{}

\begin{abstract}
Composition principles for reduced moduli are extended to the case of domains in the $n$-dimensional Euclidean space, $n>2$. As a consequence analogues of extremal decomposition theorems of Kufarev, Dubinin and Kirillova in the planer case are obtained.
\end{abstract}

{\bf Subject Classification:} 31B99  \\

{\bf Keywords:} reduced modulus, Robin function, Neumann function, nonoverlapping domain, extremal decomposition problem.

\section{Introduction and notations}

Extremal decomposition problems have a rich history and go back to
M.A.~Lavrentiev's inequality for the product of conformal radii of
non-overlapping domains. There exist two methods of their study: the
extremal-metric method and the capacitive method. The first one has
been systematically developed in papers by G.V. Kuz'mina, E.G. Emel'yanov,
A.Yu.~Solynin, A.~Vasil'ev, and Ch.~Pommerenke \cite{zbMATH01166337,zbMATH01357277,zbMATH03983566,Pom}.
The second approach is developed mainly in works of V.N.~Dubinin and
his students \cite{zbMATH01637609,zbMATH01637610,zbMATH02213855,zbMATH05659053}.
In particular, a series of well-known results about extremal decomposition
follows one way from composition principles for generalized reduced
moduli (see \cite[p. 56]{zbMATH06313560} and \cite{zbMATH06009104}). In the present paper we extend the mentioned composition
principles to the case of spatial domains. As a consequence we get
theorem about extremal decomposition for the harmonic radius \cite{zbMATH03279471}
obtained earlier in \cite{zbMATH01637610}.

Throughout the paper, $\mathbb{R}^{n}$ denotes the $n$-dimensional
Euclidean space consisting of points $x=(x_1,\dots,\,x_n)$, $n\geq3$,
and $|x|=\sqrt{x_1^2+\dots+x_n^2}$ is the length of a
vector $x\in\mathbb{R}^{n}$. We introduce the following notations:

$B(a,r)=\{x\in\mathbb{R}^{n}:\,|a-x|<r\},$

$S(a,r)=\{x\in\mathbb{R}^{n}:\,|a-x|=r\}$, ~~~$a\in\mathbb{R}^{n}$;

$\omega_{n-1}=2\pi^{n/2}/\Gamma(n/2)$ is the area of the unit sphere
$S(0,1)$;

$\lambda_{n}=((n-2)\omega_{n-1})^{-1}$.

$D$ is a bounded domain in $\mathbb{R}^{n}$, $\Gamma$ is a closed
subset of $\partial D$. The pair $(D,\Gamma)$ is admissible if there
exists the Robin function, $g_{\Gamma}(z,z_{0},D)$ harmonic in $D\setminus\{z_{0}\}$, continuous in  $\overline{D}\setminus\{z_{0}\}$ and
\begin{align}
\frac{\partial g_{\Gamma}}{\partial n} & =0\,\,\text{on}\,\,(\partial D)\setminus\Gamma,\label{eq:norm_der_rob}\\
g_{\Gamma} & =0\,\,\text{on\,\,\ensuremath{\Gamma}},\label{eq:rob}
\end{align}
and in a neighborhood of $z_{0}$ there is an expansion
\begin{equation}
g_{\Gamma}(z,z_{0},D)=\lambda_{n}\left(|z-z_{0}|^{2-n}-r(D,z_{0},\Gamma)^{2-n}+o(1)\right),\,\,z\rightarrow z_{0},\label{eq:robz0}
\end{equation}
where $\partial/\partial n$ means the inward normal derivative on the boundary. In what follows all such pairs are assumed to be admissible.

In the case $\Gamma=\emptyset$ we change the condition (\ref{eq:norm_der_rob})
by the condition
\[
\frac{\partial g_{\Gamma}}{\partial n}=\frac{1}{\mu_{n-1}(\partial D)}\,\,\text{on}\,\,\partial D,
\]
where $\mu_{n-1}(\partial D)$ is the area of boundary.

By analogy with the definition of the Robin radius for plain domains from the paper \cite{zbMATH05659053} we will call the constant $r(D,z_0,\Gamma)$ the Robin radius of the domain $D$ and the set $\Gamma$. Note that in the case of $\Gamma=\partial D$ we get the harmonic radius \cite{zbMATH03279471,Lev,zbMATH01637610}.

Let $\Delta=\{\delta_{k}\}_{k}^{m}$ be a collection of real numbers
and $Z=\{z_{k}\}_{k=1}^{m}$ be points of the domain $D$. For $\Gamma=\emptyset$
we additionally require
\[
\sum_{k=1}^{m}\delta_{k}=0.
\]
Define the potential function for the domain $D$, the set $\Gamma$,
the collection of points $Z$, and numbers~$\Delta$:
\[
u(z)=u(z;Z,D,\Gamma,\Delta)=\sum_{k=1}^{m}\delta_{k}g_{\Gamma}(z,z_{k},D).
\]

Note that for $\Gamma =\emptyset$  the function $g_{\Gamma}\left(z,z_k,D\right)$ is defined up to an additive constant. Nevertheless, the function $u(z)$ is defined uniquely and characterized by the condition
$$
\frac{\partial u}{\partial n}=0 \ \ \text{on} \ \ \partial D.
$$

It is clear from the definition of the potential function that in a neighborhood of $z_{k}$ we have
\[
u(z)=\delta_{k}\lambda_{n}|z-z_{k}|^{2-n}+a_{k}+o(1), \,\,k=1,\dots,m,
\]
where
\[
a_{k}=-\delta_{k}\lambda_{n}r(D,z_{k},\Gamma)^{2-n}+\sum_{\substack{l=1\\
l\ne k
}
}^{m}\delta_{l}g_{\Gamma}(z_{l},z_{k},D).
\]

Now if we introduce the following notation
\[
g_{\Gamma}(z_{k},z_{k},D)=-\lambda_{n}r(D,z_{k},\Gamma)^{2-n},
\]
then the constant in the expansion of the potential function in a
neighbourhood of $z_{k}$ is
\begin{equation}
a_{k}=\sum_{l=1}^{m}\delta_{l}g_{\Gamma}(z_{l},z_{k},D).\label{eq:ak}
\end{equation}

A function $v(z)$ is admissible for $D$, $Z$, $\Delta$, and $\Gamma$
if $v(z)\in\text{Lip}$ in a neighbourhood of each point of $D$ except
maybe finitely many such points, continuous in $\overline{D} \setminus \bigcup_{k=1}^{m}\{z_{k}\}$, $v(z)=0$ on $\Gamma$, and in neighbourhood
of $z_{k}$ there is an expansion
\begin{equation}
v(z)=\delta_{k}\lambda_{n}|z-z_{k}|^{2-n}+b_{k}+o(1),\,\,z\rightarrow z_{k.}\label{admiss}
\end{equation}

The Dirichlet  integral is the following
\[
I(f,D)=\int_{D}\left|\nabla f\right|^{2}d\mu,
\]
where $d\mu=dx_{1}\dots dx_{n}$.

\section{Main results}
\begin{Lemma}\label{2.1}
The asymptotic formula
$$
I(u,D_{r})=\left(\sum_{k=1}^{m}\delta_{k}^{2}\right)\lambda_{n}r^{2-n}+\sum_{k=1}^{m}\delta_{k}a_{k}+o(1),\,\,r\rightarrow0,
$$
is true, where $u$ is the potential function and $a_k$, $k=1,\dots, m$ are defined in~$(\ref{eq:ak})$ and $D_r=D\setminus \bigcup_{k=1}^{m}B(z_k,r)$.
\end{Lemma}

\textbf{Proof.}
The Green's identity
\[
\int_{V}|\nabla u|^{2}d\mu=-\int_{\partial V}u\frac{\partial u}{\partial n}ds
\]
gives
\begin{equation}
I\left(u,D_{r}\right)=-\int_{\partial D_{r}}u\frac{\partial u}{\partial n}ds=-\sum_{k=1}^{m}\int_{S(z_{k},r)}u\frac{\partial u}{\partial n}ds.\label{eq:iubr-1}
\end{equation}
The second equality in (\ref{eq:iubr-1}) holds because $u\dfrac{\partial u}{\partial n}=0$ on $\partial D$.
Note that
\[
u=\delta_{k}\lambda_{n}r^{2-n}+a_{k}+o(1),\,\,z\rightarrow z_{k}
\]
in a neighbourhood of $z_{k}$.

We calculate the integral $\int_{S(z_k,r)}u\frac{\partial u}{\partial n} ds.$
Let $u(z)=h(z)+g(z)$, where $h(z)=\lambda_n \delta_k |z-z_k|^{2-n}$ and $g(z)$ is a harmonic function. Note that $g(z_k)=a_k$. For $|z-z_k|=r$ we have the following correlations
$$
r^{n-1}u\frac{\partial u}{\partial n} =r^{n-1} \left(h\frac{\partial h}{\partial n}+h \frac{\partial g}{\partial n}+g\frac{\partial h}{\partial n}+g\frac{\partial g}{\partial n}\right)
$$
$$
=(2-n)\lambda_n^2 \delta_k^2 r^{2-n} +r\lambda_n \delta_k \frac{\partial g}{\partial n}+(2-n)g\lambda_n\delta_k+g\frac{\partial g}{\partial n}r^{n-1}
$$
$$
=-\frac{\lambda_n \delta_k^2}{\omega_{n-1}}r^{2-n}-\frac{g(z_k)\delta_k}{\omega_{n-1}}+o(1), r\rightarrow 0.
$$
Therefore
$$
\int_{S(z_k,r)}u\frac{\partial u}{\partial n}ds =\int_{S(0,1)}u\frac{\partial u}{\partial n} r^{n-1} ds
$$
$$
=-\lambda_n\delta_k^2 r^{2-n}-\delta_k a_k +o(1), \ \ r\rightarrow 0.
$$
Substituting it in (\ref{eq:iubr-1}) we get the lemma. $\square$

\begin{Lemma}\label{2.2}
For an admissible function $v$ and the potential function $u$ we
have
\[
I(v-u,D_{r})=I(v,D_{r})-I(u,D_{r})-2\sum_{k=1}^{m}\delta_{k}(b_{k}-a_{k})+o(1),\,\,r\rightarrow0.
\]

\end{Lemma}

\textbf{Proof.}
One may observe that
\begin{multline*}
I(v-u,D_{r})=\int_{D_{r}}\left(|\nabla v|^{2}+|\nabla u|^{2}-2\nabla u\nabla v\right)d\mu\\
=\int_{D_{r}}\left(|\nabla v|^{2}-|\nabla u|^{2}\right)d\mu+2\int_{\partial D_{r}}(v-u)\frac{\partial u}{\partial n}ds\\
=I(v,D_{r})-I(u,D_{r})+2\sum_{k=1}^{m}\int_{S(z_k,r)}(v-u)\frac{\partial u}{\partial n}ds\\
=I(v,D_{r})-I(u,D_{r})-2\sum_{k=1}^{m}\delta_{k}(b_{k}-a_{k})+o(1),\,\,r\rightarrow0.
\end{multline*}
Here we calculated the integral $\int_{S(z_k,r)}(v-u)\frac{\partial u}{\partial n}$ a similar way as in the proof of lemma \ref{2.1} and used the Green's identity
\[
\int_{D_{r}}(\nabla u\cdot \nabla v)d\mu=-\int_{\partial D_{r}}v\frac{\partial u}{\partial n}ds,
\]
where $n$ is the inner normal vector.
$\square$

The quantity
\[
\sum_{k=1}^{m}\delta_{k}a_{k}=\sum_{k=1}^{m}\sum_{l=1}^{m}\delta_{k}\delta_{l}g_{\Gamma}(z_{l},z_{k},D)
\]
we call the reduced modulus and denote it by $M(D,\Gamma,Z,\Delta).$
According to lemma \ref{2.1}
\[
M(D,\Gamma,Z,\Delta)=\lim_{r\rightarrow0}\left(I(u,D_{r})-\left(\sum_{k=1}^{n}\delta_{k}^{2}\right)\lambda_{n}r^{2-n}\right).
\]

\begin{Theorem}\label{2.3}
Let sets $D$, $\Gamma$, collections $Z=\{z_{k}\}_{k=1}^{m}$,
$\Delta=\{\delta_{k}\}_{k=1}^{m}$, be as in the definition of the
reduced modulus $M=M(D,\Gamma,$ $Z,\Delta),$ $u(z)$ be the potential
function for $D$, $\Gamma$, $Z$, $\Delta$, and let $D_{i}\subset D$~
be pairwise non-overlapping subdomains of $D$, $\Gamma_{i}$, $Z_{i}=\{z_{ij}\}_{j=1}^{n_{i}}$,
$\Delta_{i}=\{\delta_{ij}\}_{j=1}^{n_{i}}$, ~be from the definition
of the reduced moduli $M_{i}=M(B_{i},\Gamma_{i},$ $Z_{i},$ $\Delta_{i}),$
$u_{i}(z)$ be the potential function for $D_{i},$ $\Gamma_{i}$,
$Z_{i},$ $\Delta_{i},$ $i=1,...,p.$ Assume that the following conditions
are fulfilled:

\begin{tabular}{l}
1)$(D\cap\partial D_{i})\subset\Gamma_{i},$ $i=1,...,p;$\tabularnewline
2) $\Gamma\subset\left(\bigcup_{i=1}^{p}\Gamma_{i}\right)\bigcup\left[\mathbb{R}^{n}\setminus\left(\bigcup_{i=1}^{p}\overline{D}_{i}\right)\right]$;\tabularnewline
3)$Z=\bigcup_{i=1}^{p}Z_{i},$ \emph{that is each point $z_{k}\in Z$
coincides with some point}\tabularnewline
$z_{ij}\in Z_{i}$ \emph{ for $k=k(i,j)$ and vice versa; }\tabularnewline
4) $\delta_{k}=\delta_{ij}$ for $k=k(i,j).$\tabularnewline
\end{tabular}

Then the inequality
\[
M\geq\sum_{i=1}^{p}M_{i}+\sum\limits _{i=1}^{p}I(u-u_{i},D_{i})\geq \sum_{i=1}^{p}M_{i}
\]
holds.\end{Theorem}

\textbf{Proof.}
Consider the function
\[
v(z)=\left\{ \begin{array}{ll}
u_{i}(z),\ \  & z\in D_{i},\\
0, & z\in D\setminus\left(\bigcup_{i=1}^{p}D_{i}\right).
\end{array}\right.
\]
The condition 1) guarantees that the function $v(z)$ is continuous
in $\overline{D}\backslash\bigcup_{k=1}^{m}\{z_{k}\}$. From the conditions
2) and 3) it follows that $v(z)=0$ for $z\in\Gamma$ and in a neighbourhood
of $z_{k},$ $k=1,...,m$, there is the expansion (\ref{admiss}).
Applying lemma \ref{2.2}, we get
\begin{equation}
I(v-u,D)=I(v,D_{r})-I(u,D_{r})-2\sum\limits _{k=1}^{p}\delta_{k}(b_{k}-a_{k})+o(1),\ r\rightarrow0,\label{eq:int}
\end{equation}
here $a_k$ and $b_k$ from (\ref{eq:ak}) and (\ref{admiss}) respectively.
By lemma \ref{2.1}
\[
I(v,D_{r})=\lambda_{n}r^{2-n}\sum\limits _{i=1}^{p}\sum\limits _{j=1}^{n_{i}}\delta_{ij}^{2}+\sum\limits _{i=1}^{p}M_{i}+o(1)=\lambda_{n}r^{2-n}\sum\limits _{k=1}^{m}\delta_{k}^{2}+\sum\limits _{i=1}^{p}M_{i}+o(1),
\]
\[
I(u,D_{r})=\lambda_{n}r^{2-n}\sum\limits _{k=1}^{m}\delta_{k}^{2}+M+o(1),\qquad r\rightarrow0,
\]
taking into account 3), we have
\[
\sum\limits _{k=1}^{m}\delta_{k}(b_{k}-a_{k})=\sum\limits _{i=1}^{p}M_{i}-M.
\]
Substituting the obtained correlations in (\ref{eq:int}), we see
that the inequality
\[
\sum\limits _{i=1}^{p}I(u-u_{i},D_{i})\leq I(v-u,D)=M-\sum\limits _{i=1}^{p}M_{i}+o(1),\qquad r\rightarrow0,
\]
is true. Theorem is proved. $\square$

\begin{Theorem}\label{2.4}
Let sets $D$, $\Gamma$, collections $Z=\{z_{k}\}_{k=1}^{m}$,
$\Delta=\{\delta_{k}\}_{k=1}^{m}$, be as in the definition of the
reduced modulus $M:=M(D,\Gamma,Z,\Delta),$ $u(z)$ be the potential
function for $D$, $\Gamma$, $Z$, $\Delta$, and let $D_{i}\subset D$,
$i=1,...,p,$~ be pairwise non-overlapping domains, $\Gamma_{i}$,
$Z_{i}=\{z_{ij}\}_{j=1}^{n_{i}}$, $\Delta_{i}=\{\delta_{ij}\}_{j=1}^{n_{i}}$,
~be from the definition of the reduced moduli $M_{i}=M(D_{i},\Gamma_{i},Z_{i},\Delta_{i}),$
$u_{i}(z)$ be the potential function for $D_{i},$ $\Gamma_{i}$,
$Z_{i},$ $\Delta_{i},$ $i=1,...,p.$ Assume that $\Gamma_{i}\subset\Gamma,\ i=1,...,p,$
$Z=\bigcup_{i=1}^{m}Z_{i},$ (that is each point $z_{k}\in Z$ coincides
with some point $z_{ij}\in Z_{i}$ for $k=k(i,j)$ and vice versa),
$\delta_{k}=\delta_{ij}.$ Then the inequality
\[
\sum_{i=1}^{p}M_{i}\geq M+\sum\limits _{i=1}^{p}I(u-u_{i},D_{i})\geq M
\]
holds.\end{Theorem}

\textbf{Proof.}
The function $u$ is admissible for $D_{i}$, $i=1,\dots,p$. Let
$b_{k}$ be constants from the expansion of the function $u$ in a
neighbourhood of $z_{k},$ $b_{ij}=b_{k}$ if
$k=k(i,j).$ Applying lemmata \ref{2.1} and  \ref{2.2} with the potential functions
$u_{k}$ for $D_{k}$ we get
\begin{multline*}
\sum_{i=1}^{p}\sum_{j=1}^{n_{i}}(\delta_{ij})^{2}r^{2-n}\lambda_{n}+\sum_{i=1}^{p}\sum_{j=1}^{n_{i}}\delta_{ij}a_{ij}+o(1)=\sum_{i=1}^{p}I(u_{i},(D_{i})_{r})=\\
=\sum_{i=1}^{p}\left(I(u,(D_{i})_{r})-2\sum_{j=1}^{n_{i}}\delta_{ij}(b_{ij}-a_{ij})-I(u-u_{i},(D_{i})_{r})\right)+o(1)\\
\leq I(u,D_{r})-\sum_{i=1}^{p}I(u-u_{i},(D_{i})_{r})-2\sum_{i=1}^{p}\sum_{j=1}^{n_{i}}\delta_{ij}(b_{ij}-a_{ij})+o(1)\\
=\sum_{i=1}^{p}\sum_{j=1}^{n_{i}}(\delta_{ij})^{2}r^{2-n}\lambda_{n}+\sum_{i=1}^{p}\sum_{j=1}^{n_{i}}\delta_{ij}b_{ij}-2\sum_{i=1}^{p}\sum_{j=1}^{n_{i}}\delta_{ij}(b_{ij}-a_{ij})\\
-\sum_{i=1}^{p}I(u-u_{i},(D_{i})_{r})+o(1),\,\,r\rightarrow0.
\end{multline*}
It implies that
\[
\sum_{i=1}^{p}\sum_{j=1}^{n_{i}}\delta_{ij}b_{ij}\leq\sum_{i=1}^{p}\sum_{j=1}^{n_{i}}\delta_{ij}a_{ij}-\sum_{i=1}^{p}I(u-u_{i},D_{i})
\]
or equivalently
\[
\sum_{i=1}^{p}I(u-u_{i},D_{i})+M(D,\Gamma,Z,\Delta)\leq\sum_{i=1}^{p}M(D_{i},\Gamma_{i},Z_{i},\Delta_{i}).
\]
Here we used the fact that the function $u-u_i$ has no singularity in $D_i$.
$\square$

Denote by $r(D_{l},x_{l})=r(D_{l},x_{l},\partial D)$ the harmonic
radius. Directly from theorem \ref{2.3} we get theorem 2 of the paper \cite{zbMATH01637610}
\begin{Corollary}
For any non-overlapping domains $D_{l}\subset\mathbb{R}^{n}$, $n\geq3$,
points $x_{l}\in D_{l}$ and real numbers $\delta_{l}$, $l=1,\dots,m$
the inequality
\[
-\sum_{l=1}^{m}\delta_{l}^{2}r\left(D_{l},x_{l}\right)^{2-n}\leq\sum_{l=1}^{m}\sum_{\substack{p=1\\
p\ne l
}
}^{m}\delta_{l}\delta_{p}|x_{l}-x_{p}|^{2-n}
\]
 holds true.\end{Corollary}

\textbf{Proof.}
The Green's function of the ball $B(0,\rho)$ is
\[
\lambda_{n}\left(|x-x_{0}|^{2-n}-\left|\frac{|x_{0}|x}{\rho}-\frac{\rho x_{0}}{|x_{0}|}\right|^{2-n}\right).
\]
 Denote by $D_{l}(\rho)$ the intersection $D_{l}\cap B(0,\rho)$.
By theorem \ref{2.3}
\[
M(\rho)\geq\sum_{l=1}^{m}M_{l}(\rho),
\]
 where $M(\rho)$ is the modulus of the ball $B(0,\rho)$, the collections
$\left\{ x_{l}\right\} _{l=1}^{m}$, $\Delta=\left\{ \delta_{l}\right\} _{l=1}^{m}$,
and $\Gamma=\partial B$,
\[
M_{l}(\rho)=-\delta_{l}^{2}r\left(D_{l}(\rho),x_{l}\right)^{2-n}\lambda_n.
\]
It is sufficient to take a limit as $\rho\rightarrow\infty$.
$\square$

Theorems \ref{2.3} and \ref{2.4} imply for $p=1$ monotonicity of the quadratic form
\[
\sum_{l=1}^{m}\sum_{p=1}^{m}\delta_{l}\delta_{p}g_{\Gamma}(z_{l},z_{p},D)
\]
under extension of a domain. Following \cite{zbMATH02213855} we will
say that a domain $\tilde{D}$ is obtained by extending a domain $D$
across a part of its boundary $\gamma\subset\partial D$ if $D\subset\tilde{D}$
and $(\partial D)\cap\tilde{D}$ lies in $\gamma$.

\begin{Corollary}
If $\tilde{D}$ is obtained by extending $D$ across $\Gamma$, $\tilde{\Gamma}\subset\left(\Gamma\bigcup(\mathbb{R}^{n}\setminus \overline{D})\right),$
then for any real numbers $\delta_{l}$ and points $z_{l}\in D$
$$
\sum_{l=1}^{m}\sum_{p=1}^{m}\delta_{l}\delta_{p}g_{\tilde{\Gamma}}\left(z_{l},z_{p},\tilde{D}\right)\geq\sum_{l=1}^{m}\sum_{p=1}^{m}\delta_{l}\delta_{p}g_{\Gamma}\left(z_{l},z_{p},D\right)+I(u-\tilde{u},D)
$$
$$
\geq\sum_{l=1}^{m}\sum_{p=1}^{m}\delta_{l}\delta_{p}g_{\Gamma}\left(z_{l},z_{p},D\right).
$$
If $\tilde{D}$ is obtained by extending $D$ across the part of $(\partial D)\setminus\Gamma$,
$\tilde{\Gamma}=\Gamma$, then
$$
\sum_{l=1}^{m}\sum_{p=1}^{m}\delta_{l}\delta_{p}g_{\tilde{\Gamma}}\left(z_{l},z_{p},\tilde{D}\right)\leq\sum_{l=1}^{m}\sum_{p=1}^{m}\delta_{l}\delta_{p}g_{\Gamma}\left(z_{l},z_{p},D\right)-I(u-\tilde{u},D)
$$
$$
\leq\sum_{l=1}^{m}\sum_{p=1}^{m}\delta_{l}\delta_{p}g_{\Gamma}\left(z_{l},z_{p},D\right),
$$
 here $u$ and $\tilde{u}$ are the potential functions for $D$,
$\Gamma$, $Z=\left\{ z_{l}\right\} _{l=1}^{m}$, $\Delta=\left\{ \delta_{l}\right\} _{l=1}^{m}$
and $\tilde{D}$, $\tilde{\Gamma}$, $Z=\left\{ z_{l}\right\} _{l=1}^{m}$,
$\Delta=\left\{ \delta_{l}\right\} _{l=1}^{m}$, respectively.
\end{Corollary}

In \cite{zbMATH05659053} the notion of the Robin radius
\[
r\left(D,z_{0},\Gamma\right)=\exp\lim_{z\rightarrow
z_{0}}\left(g_{D}\left(z,z_{0},\Gamma\right)+\log|z-z_{0}|\right)
\]
 was introduced. This quantity generalized the notion of the conformal
radius. An analogue of Kufarev's theorem (see \cite{Kuf}) for
non-overlapping domains $D_{1}$, $D_{2}$ lying in the unit disk
$U$ under the condition $\left(\left(\partial D_{k}\right)\cap
U\right)\subset\Gamma_{k}\subset\partial D_{k}$, $a_{k}\in D_{k}$,
$k=1,2$ is the inequality
\[
r\left(D_{1},a_{1},\Gamma_{1}\right)r\left(D_{2},a_{2},\Gamma_{2}\right)
\leq|a_{2}-a_{1}|^{2}\left[1-\left|\frac{a_{2}-a_{1}}{1-\overline{a_{1}}a_{2}}\right|^2\right]^{-1}.
\]
By setting in theorem \ref{2.3} $p=2$, $\Gamma=\emptyset$, we obtain in
$\mathbb{R}^{n}$ the following inequality.
\begin{Corollary}
Let $D_{1}$ and $D_{2}$ be non-overlapping and lie in the ball
$U=B(0,1)$, $a_{k}\in D_{k}$, $\left(\partial D_{k}\right)\cap
U\subset\Gamma_{k}\subset\partial D_{k}$, $k=1,2$. Then
\begin{equation}
-\lambda_{n}r\left(D_{1},a_{1},\Gamma_{1}\right)^{2-n}-\lambda_{n}r\left(D_{2},a_{2},\Gamma_{2}\right)^{2-n}\leq
M\left(U,\emptyset,\left\{ a_{1},a_{2}\right\} ,\left\{
1,-1\right\} \right).\label{eq:cor3}
\end{equation}

\end{Corollary}

To calculate $M\left(U,\emptyset,\left\{ a_{1},a_{2}\right\}
,\left\{ 1,-1\right\} \right)$ we need to know the Neumann
function of the unit ball. Note that it is a quite comlicated
problem in $\mathbb{R}^{n}$. In particular, for $n=3$ (see
\cite{Sad})
\[
g_{\emptyset}\left(x,y, U\right)=
\frac{1}{4\pi}\left(\frac{1}{|x-y|}+\frac{|y|}{|x|y|^{2}-y|}-\log\left|1-(x,y)+\frac{\left|x|y|^{2}-y\right|}{\left|y\right|}\right|\right).
\]
 In \cite{Sad} there is an analytic view of $g_{\emptyset}(D,x,y)$ for $n=4,5$.
So, for $n=3$ the inequality (\ref{eq:cor3}) has the following
form
\begin{multline*}
-r\left(D_{1},a_{1},\Gamma_{1}\right)^{-1}-r\left(D_{2},a_{2},\Gamma_{2}\right)^{-1}\leq-\frac{2}{|a_{1}-a_{2}|}-
\frac{2|a_{2}|}{|a_{1}|a_{2}|^2-a_{2}|}\\
+2\log\left|1-(a_{1},a_{2})+\frac{\left|a_{1}|a_{2}|^{2}-a_{2}\right|}{\left|a_{2}\right|}\right|
+\frac{1}{1-|a_{1}|^{2}}
+\frac{1}{1-|a_{2}|^{2}}\\
-\log(4(1-|a_1|^2)(1-|a_2|^2)).
\end{multline*}

{\bf Acknowledgements.}
This work was supported by the Russian Science Foundation under grant
14-11-00022.


\begin{thebibliography}{99}

\bibitem{zbMATH06313560}{V.N. Dubinin, \em Condenser capacities and symmetrization in geometric function theory,} Birkh\"auser/Springer, Basel, 2014.

\bibitem{zbMATH02213855} {V.N. Dubinin and N.V. \ Eyrikh, Some applications of generalized condensers in the theory of analytic functions, \em Journal of Mathematical Sciences,} {\bf 133(6)} (2006), 1634-1647.

\bibitem{zbMATH05659053} {V.N. Dubinin and D.A. Kirillova, On extremal decomposition problems, \em Journal of Mathematical Sciences,} {\bf 157(4)} (2009), 573-583.

\bibitem{zbMATH01637609} {V.N. Dubinin and L.V. Kovalev, The reduced module of the complex sphere, \em Journal of Mathematical Sciences,} {\bf 105(4)} (2001), 2165-2179.

\bibitem{zbMATH01637610} {V.N. Dubinin and E.G. Prilepkina, Extremal decomposition of spatial domains, \em Journal of Mathematical Sciences,} {\bf 105(4)} (2001), 2180-2189.

\bibitem{zbMATH03983566} {E.G. Emel'yanov, On problems on the extremal decomposition, \em Zapiski Nauchnykh Seminarov Leningradskogo Otdeleniya Matematicheskogo Instituta Imeni V. A. Steklova,} {\bf 154} (1986), 76-89.

\bibitem{zbMATH03279471} {J. Hersch, Transplantation harmonique, transplantation par modules, et th\'eoremes isoperimetriques, \em Commentarii Mathematici Helvetici,} {\bf 44} (1969), 354-366.

\bibitem{Kuf} {P.P. Kufarev, On the question of conformal mappings of complementary domains, \em Dokl. Akad. Nauk SSSR,} {\bf 73} (1950), 881-884.

\bibitem{zbMATH01166337} {G.V. Kuz'mina, Methods of geometric function theory. I, II, \em St. Petersburg Mathematical Journal,} {\bf 9(3)}  (1998), 455-507; {\bf 9(5)} (1998), 889-930.


\bibitem{Lev} {B.E. Levitskii, The reduced p-module and the inner p-harmonic radius, \em Dokl. Akad. Nauk SSSR,} {\bf 316} (1991), 812-815.

\bibitem{Pom} {Ch. Pommerenke and A. Vasil'ev, Angular
derivatives of bounded univalent functions and extremal partions
of the unit disk, \em Pacific J. Math.,} {\bf 206(2)} (2002), 425-450.

\bibitem{zbMATH06009104} {E.G. Prilepkina, On composition principles for reduced moduli, \em Siberian Mathematical Journal,} {\bf 52(6)} (2011), 1079-1091.

\bibitem{Sad} {M.A. Sadybekov and B.T. Torebek and B.Kh. Turmetov, Representation of Green's function of the Neumann problem for a multi-dimensional ball, \em Complex Variables and Elliptic Equations,} http://dx.doi.org/10.1080/17476933.2015.1064402

\bibitem{zbMATH01357277} {A.Yu. Solynin, Modules and extremal metric problems, \em St. Petersburg Mathematical Journal,} {\bf 11(1)} (2000), 1–65.




\end{thebibliography}
\end{document}